\documentclass[leqno,12pt]{article}  
\usepackage{amsmath}
\usepackage{amssymb}

\usepackage[german,english]{babel}
\usepackage{amsthm}

\title{Equivariant crystalline cohomology and base change}

\author{\textsc{Elmar Grosse-Kl\"onne}}
\date{}

\theoremstyle{plain} 
\newtheorem{satz}{Theorem}[section]  
\newtheorem{lem}[satz]{Lemma}  
 
\newcommand{\spec}{\mbox{\rm Spec}}  

\newcommand{\quot}{\mbox{\rm Quot}}  
\newcommand{\spf}{\mbox{\rm Spf}}  

\newcommand{\bi}{\mbox{\rm im}}  

\theoremstyle{remark}

\theoremstyle{definition}

\begin{document}
\maketitle
\footnote[0]{}

\begin{abstract} Given a perfect field $k$ of characteristic $p>0$, a smooth proper $k$-scheme $Y$, a crystal $E$ on $Y$ relative to $W(k)$ and a finite group $G$ acting on $Y$ and $E$, we show that, viewed as virtual $k[G]$-module, the reduction modulo $p$ of the crystalline cohomology of $E$ is the de Rham cohomology of $E$ modulo $p$. On the way we prove a base change theorem for the virtual $G$-representions associated with $G$-equivariant objects in the derived category of $W(k)$-modules.
\end{abstract}

%


\section{The Theorem}

Let $k$ be a perfect field of characteristic $p>0$, let $W$ denote its ring of Witt vectors, let $K=\quot(W)$. Let $Y$ be a proper and smooth $k$-scheme and suppose that the finite group $G$ acts (from the right) on $Y$. Let $E$ be a locally free, finitely generated crystal of ${\mathcal O}_{Y/W}$-modules and suppose that for each $g\in G$ we are given an isomorphism of crystals $\tau_g:E\to g^*E$ (where $g^*E$ denotes the pull back of $E$ via $g:Y\to Y$) such that $g_2^*(\tau_{g_1})\circ\tau_{g_2}=\tau_{g_2g_1}$ (equality as maps $E\to (g_2g_1)^*E=g_2^*g_1^*E$) for any two $g_1, g_2\in G$. For $s\in{\mathbb Z}$ let $H_{crys}^s(Y/W,E)$ denote the $s$-th crystalline cohomology group (relative to $\spf(W)$) of the crystal $E$, a finitely generated $W$-module which is zero if $s\notin[0,2\dim(Y)]$ (see \cite{bo}). On the other hand, the reduction modulo $p$ of the crystal $E$ is equivalent with a locally free ${\mathcal O}_Y$-module $E_k$ with connection $E_k\to E_k\otimes_{{\mathcal O}_Y}\Omega^{1}_Y$; here $\Omega^{1}_Y$ denotes the ${\mathcal O}_Y$-module of differentials of $Y/k$. Let $\Omega^{\bullet}_Y\otimes E_k$ denote the corresponding de Rham complex. The cohomology group $H^s(Y,\Omega^{\bullet}_Y\otimes E_k)$ is a finite dimensional $k$-vector space which is zero if $s\notin[0,2\dim(Y)]$. The isomorphisms $\tau_g$ for $g\in G$ provide each $H_{crys}^s(Y/W,E)$, each $H^s(Y,\Omega^{\bullet}_Y\otimes E_k)$ and each $H^s(Y,\Omega^{t}_Y\otimes E_k)$ with an action of $G$ (from the left). By definition, the reduction modulo $p$ of the $K[G]$-module $H^s_{crys}(Y/W,E)\otimes_WK$ is the ${k}[G]$-module obtained by reducing modulo $p$ the $G$-stable $W$-lattice $H^s_{crys}(Y/W,E)/({\rm torsion})$ in $H^s_{crys}(Y/W,E)\otimes_WK$. 

\begin{satz}\label{brauer} For any $j$, the following three virtual ${k}[G]$-modules are the same:\\ (i) the reduction modulo $p$ of the virtual $K[G]$-module $\sum_s(-1)^sH^s_{crys}(Y/W,E)\otimes_WK$\\(ii) $\sum_s(-1)^sH^s(Y,\Omega^{\bullet}_Y\otimes E_k)$\\(iii) $\sum_{s,t}(-1)^{s+t}H^s(Y,\Omega^{t}_Y\otimes E_k)$.
\end{satz}

An obvious variant of Theorem \ref{brauer} holds in logarithmic crystalline cohomology, for crystals $E$ on the logarithmic crystalline site of $Y/W$ with respect to a log structure defined by a normal crossings divisor on $Y$. 
Similarly, the proof which we give below also shows the analog of Theorem \ref{brauer} for the $\ell$-adic cohomology ($\ell\ne p$) of constructible $\ell$-adic sheaves on $Y$, even if $Y/k$ is not proper. Of course, the result in the $\ell$-adic case (even for non-poper $Y/k$) is well known; it has been used for investigating the reduction modulo $\ell$ of the Deligne-Lusztig characters of groups $G={\mathbb G}({\mathbb F})$, where ${\mathbb G}$ is a reductive group over a finite field ${\mathbb F}$ of characteristic $p$. In \cite{dellus} we use the variant of Theorem \ref{brauer} in logarithmic crystalline cohomology to show that these Deligne-Lusztig characters, usually defined via $\ell$-adic cohomology of certain ${\mathbb F}$-varieties which are non-proper in general, can also be expressed through the log crystalline cohomology of suitable log crystals on suitable proper and smooth ${\mathbb F}$-varieties with a normal crossings divisor. Unfortunately, the (more geometric) proof of the $\ell$-adic analog of Theorem \ref{brauer} (due to Deligne and Lusztig, see for example \cite{caen} Lemma 12.4 and A3.15) breaks down for crystalline cohomology. On the other hand, our proof of Theorem \ref{brauer} contains a result (Theorem \ref{derthm}) on $G$-actions on strictly perfect complexes in the derived category which should be of independent interest.\\

\section{The Proof}

{\sc Proof of Theorem \ref{brauer}:} (ii)=(iii) is clear. By \cite{bo} we know that the total crystalline cohomology ${\mathbb R}\Gamma_{crys}(Y/W,E)$, as an object in the derived category $D(W)$ of the category of $W$-modules, is represented by a complex of $W$-modules of finite tor-dimension and with finitely generated cohomology; by functoriality, $G$ acts on ${\mathbb R}\Gamma_{crys}(Y/W,E)$. Also from \cite{bo} we know that the total crystalline cohomology commutes with base change, i.e. that ${\mathbb R}\Gamma_{crys}(Y/W,E)\otimes_W^{\mathbb L}k$ is the total crystalline cohomology of the reduction modulo $p$ of $E$ (as a crystal relative to $\spec(k)$). But the latter is known (see \cite{bo} Corollary 7.4) to be the de Rham cohomology of $E_k$, i.e. its $s$-th cohomology group is $H^s(Y,\Omega^{\bullet}_Y\otimes E_k)$. Hence (i)=(ii) follows from Theorem \ref{derthm} below.\hfill$\Box$\\

Let $A$ be a complete discrete valuation ring with perfect residue field $k$ of characteristic $p>0$ and fraction field $K$ of characteristic $0$. Let $L^{\bullet}$ be a complex of $A$-modules of finite tor-dimension and with finitely generated cohomology; by \cite{bo} Lemma 7.15 this is equivalent with saying that $L^{\bullet}$ is quasiisomorphic to a strictly perfect complex, i.e. a bounded complex of finitely generated projective $A$-modules. Suppose the finite group $G$ acts on $L^{\bullet}$ when $L^{\bullet}$ is viewed as an object in the derived category $D(A)$ of the category of $A$-modules. Then each cohomology group $H^i(L^{\bullet}\otimes_A K)=H^i(L^{\bullet})\otimes_A K$ (resp. each cohomology group $H^i(L^{\bullet}\otimes_A^{\mathbb L}k)$) becomes a representation of $G$ on a finite dimensional $K$-vector space (resp. $k$-vector space). 

\begin{satz}\label{derthm} The virtual $k[G]$-module $\sum_i(-1)^iH^i(L^{\bullet}\otimes_A^{\mathbb L}k)$ is the reduction (modulo the maximal ideal of $A$) of the virtual $K[G]$-module $\sum_i(-1)^iH^i(L^{\bullet})\otimes_AK$. Equivalently, the restriction of the character of $\sum_i(-1)^iH^i(L^{\bullet})\otimes_AK$ to the subset of $p$-regular elements of $G$ is the Brauer character of $\sum_i(-1)^iH^i(L^{\bullet}\otimes_A^{\mathbb L}k)$.
\end{satz}

We say that the automorphism $\gamma$ of the finitely generated $A$-module $M$ is {\it prime to $p$} if and only if the following holds. For any finite extension $A'\supset A$ with a discrete valuation ring $A'$ and for any two $\gamma\otimes_AA'$-stable submodules $N, N'$ of $M\otimes_AA'$ with $N'\subset N$ and such that $N/N'$ is a cyclic $A'$-module, the endomorphism which $\gamma\otimes_AA'$ induces on $N/N'$ is of finite order prime to $p$. 

\begin{lem}\label{criffo} Let $\gamma$ be an automorphism of the finitely generated $A$-module $M$.\\(a) If $M$ is free then $\gamma$ is prime to $p$ if and only if the roots of the characteristic polynomial of $\gamma$ are roots of unity of order prime to $p$. In particular, $\gamma|_N:N\to N$ is prime to $p$ for each submodule $N$ of $M$ with $\gamma(N)=N$.\\(b) Let $M_1\subset M$ be a submodule with $\gamma(M_1)=M_1$ and such that $M_2=M/M_1$ is free. Let $\gamma_1$, resp. $\gamma_2$, be the induced automorphism of $M_1$, resp. of $M_2$. If $\gamma_1$ and $\gamma_2$ are prime to $p$, then $\gamma$ is prime to $p$.
\end{lem}

{\sc Proof:} Statement (a) is clear. (b) Let $N'\subset N\subset M\otimes_AA'$ be as in the definition. If $N\subset M_1\otimes_AA'$ the hypothesis on $\gamma_1$ applies. Otherwise, since $M_2\otimes_AA'$ is free over $A'$ and $N/N'$ is cyclic, $N/N'$ maps injectively to $M_2\otimes_AA'$ and the hypothesis on $\gamma_2$ applies.\hfill$\Box$\\ 

{\sc Proof of Theorem \ref{derthm}:} The problem is of course that the $H^i(L^{\bullet})$ may have torsion, i.e. $H^i(L^{\bullet})\otimes_Ak\ne H^i(L^{\bullet}\otimes_A^{\mathbb L}k)$ in general. Similarly, the task would be easy if we knew that there is a strictly perfect complex $K^{\bullet}$ quasiisomorphic to $L^{\bullet}$ such that the action of $G$ on $L^{\bullet}$ in $D(A)$ is given by the action of $G$ on $K^{\bullet}$ by true morphisms of complexes (not just by morphisms in $D(A)$).
 We introduce some notations. For an automorphism $\gamma:L^{\bullet}\to L^{\bullet}$ in $D(A)$ let $\epsilon_1^i,\ldots,\epsilon_{n(i)}^i$ (with $n(i)=\dim_kH^i(L^{\bullet}\otimes_A^{\mathbb L}k)$) denote the roots of the characteristic polynomial of $\gamma$ acting on $H^i(L^{\bullet}\otimes_A^{\mathbb L}k)$ and let $\widetilde{\epsilon}_1^i,\ldots,\widetilde{\epsilon}_{n(i)}^i$ denote their Teichm\"uller liftings. On the other hand let $\xi_1^i,\ldots\xi_{n'(i)}^i$ (with $n'(i)=\dim_KH^i(L^{\bullet})\otimes_AK$) denote the roots of the characteristic polynomial of $\gamma$ acting on $H^i(L^{\bullet})\otimes_AK$. Then let $$Br(\gamma,H^{\heartsuit}(L^{\bullet}\otimes_A^{\mathbb L}k))=\sum_i(-1)^i\sum_{j=1}^{n(i)}\widetilde{\epsilon}_j^i,$$$$Tr(\gamma,H^{\heartsuit}(L^{\bullet})\otimes_AK)=\sum_i(-1)^i\sum_{j=1}^{n'(i)}\xi_j^i.$$What we must show is that for all $p$-regular elements $g\in G$ (those whose order in $G$ is not divisible by $p$), if $\gamma:L^{\bullet}\to L^{\bullet}$ denotes the corresponding automorphism of $L^{\bullet}$ in $D(A)$, then$$Br(\gamma,H^{\heartsuit}(L^{\bullet}\otimes_A^{\mathbb L}k))=Tr(\gamma,H^{\heartsuit}(L^{\bullet})\otimes_AK).$$Clearly it is enough to show the following statement. For any strictly perfect complex $L^{\bullet}$ of $A$-modules (not necessarily endowed with a $G$-action in $D(A)$) and for any automorphism $\gamma:L^{\bullet}\to L^{\bullet}$ in $D(A)$ which on the cohomology modules induces automorphisms prime to $p$ we have$$Br(\gamma,H^{\heartsuit}(L^{\bullet}\otimes_A^{\mathbb L}k))=Tr(\gamma,H^{\heartsuit}(L^{\bullet})\otimes_AK).$$We use induction on the minimal $m\in{\mathbb Z}_{\ge0}$ with the following property: after a suitable degree shift we have $L^i=0$ for all $i\notin[0,m]$. For $m=0$ the statement is clear from Lemma \ref{criffo} (a). Now let $m\ge1$; shifting degrees we may assume $L^i=0$ for all $i\notin[0,m]$. Let $d^m:L^{m-1}\to L^m$ denote the differential. Choose a sub-$k$-vector space $N_k^{m-1}$ of $L^{m-1}\otimes k$ which under $d^m\otimes k$ maps isomorphically to the kernel of $L^m\otimes k\to H^m(L^{\bullet}\otimes k)=H^m(L^{\bullet})\otimes k$. Then $N_k^{m-1}=N^{m-1}\otimes k$ for a direct summand $N^{m-1}$ of $L^{m-1}$. By construction, $N^{m-1}$ maps isomorphically to its image $N^m$ in $L^m$. Thus, setting $N^i=0$ if $i\notin\{m-1,m\}$, the subcomplex $N^{\bullet}$ of $L^{\bullet}$ is acyclic. Dividing it out we may therefore assume $L^m\otimes k=H^m(L^{\bullet}\otimes k)$.
Since the functor $K^-({\rm proj-}A)\to D(A)$ from the homotopy category of complexes of projective $A$-modules bounded above to $D(A)$ is fully faithful, the action of $\gamma$ on $L^{\bullet}$ in $D(A)$ is in fact represented by a true morphism of complexes $\gamma^{\bullet}:L^{\bullet}\to L^{\bullet}$. Base changing to a finite extension of $A$ by a discrete valuation ring (this does not affect the numbers $Br$ and $Tr$) we may suppose that the characteristic polynomial of $\gamma^m:L^m\to L^m$ splits in $A$ (we remark that $\gamma^m$ is bijective: this follows from $L^m\otimes k=H^m(L^{\bullet}\otimes k)$ and the fact that $\gamma$ acts bijectively on $H^m(L^{\bullet}\otimes k)$). We therefore find a $\gamma^m$-stable filtration\begin{gather}(0)=F^0\subset F^1\subset\ldots\subset F^s=L^m\quad\quad(s={\rm {rk}}(L^m))\label{fil}\end{gather}such that $G^e=F^e/F^{e-1}$ is free of rank one, for any $1\le e\le s$. The cyclic $A$-module$$\frac{F^e}{(F^e\cap\bi(d^m))+F^{e-1}}$$is a $\gamma^m$-stable subquotient of $H^m(L^{\bullet})$ (it is non-zero because of $L^m\otimes k=H^m(L^{\bullet}\otimes k)$), hence $\gamma^m$ acts on it by multiplication with a root of unity of order prime to $p$. Let $\xi_e\in A^{\times}$ denote its Teichm\"uller lifting. Choose $\ell_e\in F^e$ which represents a basis element of $G^e$; then $\ell_1,\ldots,\ell_s$ is a basis of $L^m$. Modulo $F^{e-1}$ the class of $\xi_e\ell_e-\gamma^m(\ell_e)\in F^e$ lies in $\bi(d^m)$. Choose a $t_e\in L^{m-1}$ with $d^m(t_e)=\xi_e\ell_e-\gamma^m(\ell_e)$ modulo $F^{e-1}$. Let $t:L^m\to L^{m-1}$ denote the $A$-linear map which sends $\ell_e$ to $t_e$, for each $1\le e\le s$. Using $t$ we see that we may modify $\gamma^{\bullet}$ within its homotopy class to achieve that the filtration (\ref{fil}) is still $\gamma^m$-stable and such that $\gamma^m$ acts on each $G^e$ by multiplication with a root of unity of prime-to-$p$-order in $A^{\times}$. Therefore we may assume that $\gamma^m:L^m\to L^m$ is prime to $p$. Let $L_1^m=L^m$ and $L_1^i=0$ for $i\ne m$. Then $L_1^{\bullet}$ is a $\gamma^{\bullet}$-stable subcomplex of $L^{\bullet}$ and since $Br(\gamma)$ and $Tr(\gamma)$ are additive in exact $\gamma^{\bullet}$-equivariant sequences of complexes it suffices to show $Br(\gamma)=Tr(\gamma)$ for the complexes $L_1^{\bullet}$ and $L^{\bullet}/L_1^{\bullet}$. Since these complexes are shorter than $L^{\bullet}$ this follows from the induction hypothesis. Indeed, the prime to $p$ hypothesis is clearly satisfied for $L_1^{\bullet}$ so it remains to show that $\gamma^{\bullet}$ induces automorphisms prime to $p$ on the cohomology modules of $L^{\bullet}/L_1^{\bullet}$. In degrees smaller than $m-1$ this is clear from the corresponding hypothesis on $L^{\bullet}$, only $H^{m-1}(L^{\bullet}/L_1^{\bullet})$ is critical. But $H^{m-1}(L^{\bullet})$ is a submodule of $H^{m-1}(L^{\bullet}/L_1^{\bullet})$ and the quotient $Q=H^{m-1}(L^{\bullet}/L_1^{\bullet})/H^{m-1}(L^{\bullet})$ maps isomorphically to a submodule of $L_1^m=L^m$. By Lemma \ref{criffo} (b) it suffices to show that $\gamma^{\bullet}$ induces automorphisms prime to $p$ on $H^{m-1}(L^{\bullet})$ and on $Q$. For $H^{m-1}(L^{\bullet})$ this holds by hypothesis, for $Q$ this follows from Lemma \ref{criffo} (a).\hfill$\Box$\\


\end{document}